# The continuity and uniqueness of the value function of the hybrid optimal control problem with reach time to a target set

Myong-Song Ho
Kim Il Sung University, Pyongyang DPRK
ryongnam5@yahoo.com

Kwang-Nam O
Middle School No.1, Pyongyang DPRK
ryongnam6@yahoo.com

Chol-Jun Hwang
Kim Il Sung University, Pyongyang DPRK
ryongnam27@yahoo.com

**Abstract:** The hybrid optimal control problem with reach time to a target set is addressed and the continuity and uniqueness of the associated value function is proved.

Hybrid systems involves interaction of different types of dynamics: continuous and discrete dynamics. The state of a continuous system is evolved by an ordinary differential equation until the trajectory hits the predefined jump sets: an autonomous jump set $A$ and a controlled jump set $C$. At each jump the trajectory is moved discontinuously to another Euclidean space by a discrete system.

We study the hybrid optimal control problem with reach time to a target set, prove the continuity of the associated value function $V$ with respect to the initial point under the assumption that $V$ is lower semicontinuous on the boundary of a target set, and also characterize it as an unique solution of a quasi-variational inequality in a viscosity sense using the dynamic programming principle.

## 1. Introduction

Hybrid control systems are control systems that involve both continuous and discrete dynamics and continuous and discrete controls. The hybrid system has been studied widely in the number of literatures and also has been applied successfully to address the problems in air traffic control [4], automotive control [6], chemical process control [7] and so on.

In [5], Branicket, Borkar, and Mitter introduced an unified model for hybrid control. They have specified priori regions $A_i$, $C_i$, $D_i$, $(i \in Z^+)$, namely, an autonomous jump set, a controlled jump set, and a jump destination set. In their work the state evolves according to an ordinary differential equation in some state space and when the trajectory hits a jump set the state is jumped to a destination set $D$ according one of the two possibilities: if the trajectory hits an autonomous jump set, then the state must jump while the trajectory hits a controlled jump set, the controller chooses to (it does not have to) move the trajectory discontinuously.

They proved a right continuity of the value function of infinite horizon optimal control problem and derived quasi-variational inequality satisfied by the value function and summarized algorithms for solving this inequalities.

In [1], Dharmatti and Ramaswamy studied the same problem as the one considered in [5] in the most general case when the autonomous jump set A is nonempty and the controlled jump set can be arbitrary, proved the local Holder continuity of the value function, and demonstrated that this value function is an unique solution of the quasi-variational inequality (QVI) in a viscosity sense when it is bounded.

To prove the continuity they estimated the distances between the trajectories starting from two closer points, through their continuous evolutions, and through their discrete jumps to show that the distance remains small for initial points sufficiently close enough. This allowed to getting a Holder exponent for the continuity of the value function.

They also considered in [2] the finite horizon hybrid control problem and proved continuity and uniqueness of the value function in the case that the cost functionals are assumed to be unbounded.

In this paper we observe the hybrid optimal control problems with reach time to a target set.

Unlike the infinite horizon control problem which is noted in [1], in our hybrid optimal control problem the evolution of the state is stopped and the payoff is computed when the trajectory hits the target set defined in some state space.

Our aim is to prove the continuity of the value function and to characterize the value function $V$ as the unique viscosity solution of the quasi-variational inequality.

The proof of continuity of the value function in this case is harder than in infinite horizon hybrid optimal control problem because of non-zero terminal cost $h$.

In order to have a continuous solution for the value function we have to assume that $V$ is lower semicontinuous on the boundary of the target set $\Gamma$.

Under this assumption we prove that the continuity at all points on $\partial \Gamma$ implies continuity every where.

Although the basic estimations are similar to those considered in [1], we shoud take care to the evaluation of the differences of the hitting times in the proof of the continuity.

In [1], they defined the first hitting time of the trajectory as $T(x) = \inf_{u} \{t > 0 \,|\, X_x(t, u) \in A\}$, ($X_x(t, u)$ is the state of the trajectory starting from an initial point $x$ and evolving with continuous control $u$), and proved that this is locally Lipschitz continuous (in Theorem 3.1) to investigate the distance between the first hitting points (Lemma 3.2).

On the other hand in the proof of Theorem 3.5 only two types of integrals are considered to sum up the total costs from the conclusion that the trajectories starting from two neighbouring points hit the autonomous jump set sequentially (expression 3.22).

When taking into consideration the problem of our own style, instead of using $T(x)$ we define the first hitting time of the trajectory starting from the initial point $x$ with control $u$ to evaluate the differences between two trajectories.

What matters in regard to us in the proof of the continuity is that the claculating of the costs in various time intervals because it doesn`t seem to be suitable to use the above relations between the hitting times (expression 3.22) directly by means of our problem.

This paper is organized as follows:

In section 2 we introduce the notations to define the hybrid control system and associated reach time problem to a target set and elucidate all assumptions used in deriving the results of the sequel.

In section 3 we prove the continuity of the value function with respect to the initial point under the assumption that the value function is lower semicontinuous on the boundary of a target set.

Section 4 deals with the dynamic programming principle for the present problem and the quasi-variational inequality that the value function must be satisfied in the viscosity sense.

In section 5 we prove the comparison between any two viscosity solutions to take up the issue o f uniqueness.

## 2. Mathematical preliminaries.

In this section we define the hybrid optimal control system and give all assumptions which are needed for the consideration of the continuity and uniqueness of the value function.

Hybrid control systems are dynamical systems that involve the interaction of different types of dynamics: continuous flows and discrete jumps. Therefore the state of the hybrid control system is described by a pair of the variables: a part of the state takes values in Euclidean space while another part takes values in a discrete set.

The evolution of the state of the continuous system is given as a following ordinary differential equation

$$\dot{y}(t) = f(y(t), q(t), u(t)) \tag{2.1}$$

$$y(0) = x, \tag{2.2}$$

Here $y(t) \in \Omega$, $\Omega = \bigcup_i \Omega_i$ with $\Omega_i$ is a closed connected subset of $R^{d_i}$, and $q(t) \in I = \{1, 2, \cdots\}$ is a discrete variable that is $q(t) = i$ whenever $y(t) \in \Omega_i$. At this time (2.1) can be written in another form

$$\dot{y}(t) = f_i(y(t), u(t)),$$

where $f_i : \Omega_i \times \mathcal{U} \to \Omega_i$ and $\mathcal{U}$ is control set

$$\mathcal{U} = \{a : [0, \infty) \to U \mid u(\cdot) : \text{measurable}, U : \text{compact metric space}\}.$$

Denote the state of the hybrid control system at time $t$ by $(y_x(t, u), q(t))$.

Here $y_x(t, u)$ is a solution of (2.1), (2.2).

In the hybrid control system the trajectory undergoes discrete jump to a destination set $D$ when it hits predefined autonomous jump set $A$ or controlled jump set $C$.

$$A = \bigcup_i A_i, \ A_i \subset \Omega_i \subset R^{d_i}$$

$$C = \bigcup_i C_i, \ C_i \subset \Omega_i \subset R^{d_i}$$

$$D = \bigcup_i D_i, \ D_i \subset \Omega_i \subset R^{d_i}$$

When the trajectory starting from the initial state $(x, q)$ hits the autonomous jump set $A_q$ then it is moved to a destination set $D_{q'}$ according to the transition map $g$. This map uses discrete controls to move the trajectory.

When the trajectory hits the controlled jump set $C_q$, every time the controller can choose either to jump or not to jump and if the controller chooses to jump then the trajectory is moved to a destination set $D_{q'}$ possibly in another space $\Omega_{q'}$ to continue the evolving under the continuous dynamics until it hits again the jump sets or a target set.

Then we get a sequence of hitting times of $A$ (we denote it by $\{\tau_i\}$) and a sequence of jumping times of $C$, where the controller chooses to make a jump (we denote it by $\{\xi_i\}$).

Denote the state that the trajectory hits $A$ at time $\tau_i$ by $(x_i, q_i) = (y(\tau_i^-, u), q(\tau_i^-))$, (in general, we take the trajectory to be left continuous) and the shifted state under the transition map $g : A \times I \times V_1 \to D$ by $(x_i', q_i') = (y(\tau_i^+, u), q(\tau_i^+))$.

Also for the controlled jump set $C$ we denote the state where the controller chooses to jump at time $\xi_i$ by $(y(\xi_i), q(\xi_i)) = (y(\xi_i^-, u), q(\xi_i^-))$, and jumped state by

$$(y(\xi_i)', q(\xi_i)') = (y(\xi_i^+, u), q(\xi_i^+)).$$

Now we give the following topology on $\Omega \times I$.

$(x_n, q_n)$ converges to $(x, q)$ if for any $\varepsilon > 0$, there exist some large number $N$ and some $i$ such that $q_n = q$, $x, x_n \in \Omega_i \subseteq R_i$ and $\|x_n - x\| < \varepsilon$ for $\forall n \geq N$. Here $\|\cdot\|$ means a Euclidean norms in $R^{d_i}$

Throughout the whole paper we deal with the sets $A$, $C$, $D$ and the map $f$, $g$ satisfying following assumptions.

($A_1$) For any $i \in I$, $\Omega_i$ is a closure of opened connected subset of $R^{d_i}$.

($A_2$) $A_i$, $C_i$, $D_i$ are closed, $\partial A_i$, $\partial C_i$ are $C^2$, and $\partial A_i \supseteq \partial \Omega_i$ for any $x_i \in D_i$.

($A_3$) The transition map $g: A \times I \times V_1 \to D$ is bounded, Lipschitz continuous with Lipschitz constant $G$, with the understanding $g = \{g_i\}$: if $q(t) = i$ then the transition map in $A_i$ is regarded as $g_i : A_i \times V_1 \to D_j$, where $V_1$ is discrete control set.

($A_4$) Vector field $f: \Omega \times I \times U \to \Omega$ is a Lipschitz continuous with Lipschitz constant $L$ in the state variable $x$, uniformly continuous in control variable $u$, and for any $i \in I$,
$$|f_i(x, u)| \leq F, \forall x \in \Omega_i, u \in U$$

($A_5$) (transversality condition) For any $i \in I$, $A_i$ is compact and there exist some constant $\xi_0 > 0$ such that $f_i(x_0, u) \cdot \eta(x_0) \leq -2\xi_0$, $\forall x_0 \in \partial A_i, u \in U$, where $\eta(x)$ is an outward normal vector to $\partial A_i$ at $x$.

We can give the same assumption to the controlled jump set $C$.

($A_6$) $\inf d(A_i, C_i) \geq \beta > 0$, $\inf d(A_i, D_i) \geq \beta > 0$, where $d(A, B)$ is a distance between the sets $A$ and $B$.

($A_7$) $U$, $V_1$ are the compact metric spaces.

As a matter of fact, all the assumptions above are the same as those in [1].

Next for some $j \in I$ we define a set $\Gamma \subset \Omega_j$, namely, a target set. This set satisfies the following assumptions.

($A_8$) $\Gamma \subset \Omega_j$ is a closed, $\partial \Gamma$ is compact, and
$$d(\Gamma, D_j) \geq \beta > 0, d(\Gamma, A_j) \geq \beta > 0, d(\Gamma, C_j) \geq \beta > 0$$

($A_9$) the transversality condition is satisfied on the target set $\Gamma$.

In this paper we observe the following hybrid optimal control problem with reach time to a target set $\Gamma$ for the control $(u(\cdot), v, \xi_i, y(\xi_i)')$.

The total discount is given by

$$J(x, q, u(\cdot), v, \xi_i, y(\xi_i)') = \int_0^{t_x(u)} K(y_x(t, u), q(t), u(t)) e^{-\lambda t} dt + \sum_{\tau_i < t_x(u)} C_a(x_i, q_i, v) e^{-\lambda \tau_i}$$
$$+ \sum_{\xi_i < t_x(u)} C_c(y(\xi_i), q(\xi_i), y(\xi_i)', q(\xi_i)') e^{-\lambda \xi_i} + e^{-\lambda t_x(u)} h(y_x(t_x(u), u)), \quad (2.3)$$

where $\lambda$ is a positive constant and
$$K: \Omega \times I \times U \to R,$$
$$C_a: A \times I \times V_1 \to R,$$
$$C_c: C \times I \times D \times I \to R,$$
$$h: \Gamma \to R,$$
are cost functions.

Unlike the problem in [1], our problem has an initial state in $\Omega \setminus \Gamma$ whose dynamics stopped and the payoff computed when the trajectory reaches $\Gamma$, therefore we define a reach time for $u \in \mathcal{U}$

$$t_x(u) = \begin{cases} +\infty & \{t \mid y_x(t, u) \in \Gamma\} = \phi \\ \min\{t \mid y_x(t, u) \in \Gamma\} & \{t \mid y_x(t, u) \in \Gamma\} \neq \phi \end{cases} \quad (2.4)$$

The value function $V$ of this problem is given as

$$V(x, q) = \inf_{\theta \in (U \times V_1 \times [0, +\infty) \times D)} J\left(x, q, u(\cdot), v, \xi_i, y(\xi_i)'\right). \tag{2.5}$$

Here the following assumptions are satisfied for $h, K, C_a, C_c$.

($A_{10}$)   $h \in C(\Gamma), h(x) \geq 0, x \in \Gamma$

($A_{11}$)   $K$ is a Lipschitz continuous function with Lipschitz constant in the state variable $x$ and uniformly continuous in the control variable $u$. Also
$$|K(x, q, u)| \leq K_0, \forall (x, q, u) \in \Omega \times I \times U$$

($A_{12}$)   $C_a, C_c$ are uniformly continuous in the variables $x, u$ and bounded below. Moreover $C_a$ is a Lipschitz continuous function with Lipschitz constant $C_1$ in the state variable $x$ and
$$|C_a(x, q, v)| \leq C_0, \forall (x, q, v) \in \Omega \times I \times V_1. \text{ Also}$$
$$C_c(x, q, y, q') < C_c(x, q, z, q') + C_c(z, q', y, q'), \forall x \in C_q, z \in D \cap C_{q'}, y \in D_{q'}$$

Now we define the first hitting time with control $u$ and singed distance function
$$t(x, q, u) = \inf\{t > 0 \mid y(t) \in A_q, y(0) = x, \dot{y}(t) = f_q(y(t), u(t))\}$$

$$d(x, q) = \begin{cases} -d(x, \partial A_q) & x \in \text{int } A_q \\ 0 & x \in \partial A_q \\ d(x, \partial A_q) & x \in \text{out} A_q \end{cases}.$$

Here $\text{int } A_q = A_q \setminus \partial A_q$ and $\text{out} A_q = \Omega_q \setminus A_q$.

*Remark1. The first hitting time we defined above is a different concept with the first hitting time defined in the literature [1]. (See section 3 of [1], 1263 p.)*

In the sequel discussion we replace $(y_x(t, u), q(t)), t(x, q, u), d(x, q)$ simply as $y_x(t, u), t(x, u), d(x)$, while the discrete jump to another space is not considered.

In next sections we are focused on the continuity of value function of a hybrid optimal control problem described above and the partial differential equation satisfied by the value function.

### 3. Continuity of the value function.

From the ordinary differential equation theory we can see the system whose dynamics is given in (2.1), and (2.2) has the following property
$$|y_x(t, u) - y_z(t, u)| < e^{Lt}|x - z| \tag{3.1}$$
$$|y_z(t, u) - y_z(\bar{t}, u)| < F|t - \bar{t}| \tag{3.2}$$

Here $F, L$ are the constants given in assumption ($A_4$) respectively.

Without loss of generality, let the initial points $x, z$ is in $\Omega_1$ and fix the control $u$.

For the proof of the value function we consider how the first hitting times with control $u$ depends on the initial points.

**Theorem1:** Assume ($A_1$)-($A_7$) to be satisfied. Let the first hitting times of trajectories starting from $x, z$ and moving with the fixed control $u$ be $\tau_1 = t(x, u), \lambda_1 = t(z, u)$ and the points where the trajectories hit $A$ be $x_1 = y_x(\tau_1, u), z_1 = y_z(\lambda_1, u), x_1, z_1 \in A_1$.

If $|x - z| < \delta_1$, for some $\delta_1 > 0$, then
$$|x_1 - z_1| < (1 + FC) e^{L(\tau_1 \vee \lambda_1)} |x - z|, \tag{3.3}$$
where $a \vee b = \max\{a, b\}$

(**Proof**) First, we show the following lemma for the evaluation of $|\tau_1 - \lambda_1|$.

**Lemma1:** there exist some $\delta^* > 0$, $C > 0$ such that for the control $u$
$$t(x, q, u) < Cd(x, q), \quad \forall x \in B(A_q, \delta^*) \setminus \text{int } A_q.$$

Here $B(A_q, \delta^*)$ is a neighborhood of $A_q$ and $d$ is a signed distance function.

(**Proof of lemma**) Let the trajectory starting from $x$ hit on $A_q$ at $x_0$.

From the assumption ($A_2$), there exist some $R > 0$ such that $d$ is $C^1$ in $B(\partial A_q, R)$.

Now we show that there exist some $r_0 < R$ such that the following inequality holds for any $0 \leq t < t(x, q, u)$ and $y_x(t, u) \in B(A_q, r_0)$:
$$f_q(y_x(t, u), u(t)) \cdot D_x d(y_x(t, u), q) < -\xi_0.$$

In fact from the transversality condition ($A_5$) it easy to know that for any $x_0 \in \partial A_q$ and $u_0 \in U$, there exist some $r_0 < R$ such that
$$f(x, u_0) \cdot D_x d(x) < -\xi_0, \quad \forall x \in B(x_0, r_0), \tag{*}$$
where $D_x d(x, q)$ means the derivatives of $d$ by means of the spatial variable $x$.

We can give the $r_0$ in above inequality regardless of $x_0$.

As a matter of fact, whenever we fix the control $u \in U$, we can choose $x_i \in \partial A_q$, $i = 1, 2, \cdots, n$ and corresponding $r_i > 0 (i = 1, 2, \cdots, n)$ such that $\bigcup_{i=1}^{n} B\left(x_i, \dfrac{r_i}{2}\right) \supset \partial A_q$ and
$$f(x, u_0) \cdot Dd(x) < -\xi_0, \quad \forall x \in B(x_i, r_i), \quad \forall i = \overline{1, n}$$
by the compactness of $\partial A_q$.

Lets $r_0 = \min\{\dfrac{r_1}{2}, \dfrac{r_2}{2}, \cdots, \dfrac{r_n}{2}\}$. If $x \in B(\partial A_q, r_0)$, then there exist some $x' \in \partial A_q$ and some $j \in \{1, 2, \cdots, n\}$ such that $x' \in B\left(x_j, \dfrac{r_j}{2}\right)$ and $|x - x'| < r_0$.

From this we can get $|x - x_j| \leq |x - x'| + |x' - x_j| \leq r_j$ which means (*) holds for every $x \in B(\partial A_q, \delta^*)$ that is
$$f(x, u_0) \cdot Dd(x) < -\xi_0, \quad \forall x \in B(\partial A_q, \delta^*), \tag{**}$$
where $\delta^* = \min\{r_0, \dfrac{\beta}{2}\}$.

If there exist some $x \in B(\partial A_q, r_0)$ and some $t'(0 \leq t' < t(x, q, u))$ such that
$$f_q(y_x(t, u), u(t)) \cdot D_x d(y_x(t, u), q) \geq -\xi_0$$
for $y_x(t', u) \in B(A_q, \delta^*)$, then it contradicts to (**) because it does not hold for $x = y_x(t', u)$, $u_0 = u(t')$.

Therefore if $y_x(t, u) \in B(A_q, \delta^*)$ for any $0 \leq t < t(x, q, u)$ then
$$f_q(y_x(t, u), u(t)) \cdot D_x d(y_x(t, u), q) < -\xi_0.$$

It means that $d(y_x(s, u), q) - d(x, q) = \int_0^s f_q(y_x(t, u), u(t)) \cdot D_x d(y_x(t, u), q) dt < -\xi_0 s$

Now if we choose $s_x = \dfrac{d(x, q)}{\xi_0}$, then

$$d(y_x(s_x, u), q) - d(x, q) < -d(x, q), \text{ and this implies } y_x(s_x, u) \in \text{int } A_q$$

Therefore for $C = \dfrac{1}{\xi_0}$, we have $t(x, u(t)) < s_x = \dfrac{d(x, q)}{\xi_0}$. □

Now we apply the lemma to get (3. 3).

Without loss of generality let $\lambda_1 > \tau_1$.

Define $\delta_1 = e^{-L\tau_1} \delta^*$ where $\delta^*$ is listed in the Lemma.

If $|x - z| < \delta_1$, then from (3. 1)

$$|y_x(\tau_1, u) - y_z(\tau_1, u)| < e^{L\tau_1} \delta^* e^{-L\tau_1} < \delta^*$$

Hence by using the lemma and $y_x(\tau_1, u) \in A$ we have

$$\lambda_1 - \tau_1 = t(y_z(\tau_1, u), u) < C d(y_z(\tau_1, u)) < C |y_x(\tau_1, u) - y_z(\tau_1, u)| < C e^{L\tau_1} |x - z|$$

If $\lambda_1 < \tau_1$ then we have $\tau_1 - \lambda_1 < C e^{L\lambda_1} |x - z|$ in a similar way, therefore

$$|\tau_1 - \lambda_1| < C e^{L(\tau_1 \vee \lambda_1)} |x - z| \tag{3. 4}$$

On the other hand

$$|x_1 - z_1| < |y_x(\tau_1, u) - y_z(\lambda_1, u)| \leq |y_x(\tau_1, u) - y_z(\tau_1, u)| + |y_z(\tau_1, u) - y_z(\lambda_1, u)|$$

$$< e^{L\tau_1}|x - z| + F|\tau_1 - \lambda_1| \leq e^{L(\tau_1 \vee \lambda_1)}(1 + FC)|x - z|.$$

In the third inequality we used (3. 1) and (3. 2) and last inequality is from (3. 4) □

Now without loss of generality let trajectories be moved from $x_1$, $z_1$ to

$x_1' = g_1(x_1, v)$, $z_1' = g_1(z_1, v)$ which belong to the initial set of $D_2 \subset \Omega_2$ and the evolving is going on $\Omega_2$ with fixed control $u$.

Denote the next hitting times of trajectories when they hit $A$ once again by $\tau_2$, $\lambda_2$.

The following theorem estimates $|\tau_2 - \lambda_2|$.

**Theorem 2:** Assume (2. 1)-(2. 7). Let the next hitting times of trajectories starting from $x$, $z$ and evolving with fixed control $u$ be $\tau_2$, $\lambda_2$.

Define $x_2 = y_{x_1'}(\tau_2 - \tau_1, u)$, $z_2 = y_{z_1'}(\lambda_2 - \lambda_1, u)$.

Then there exist some $\delta_2 > 0$ such that if $|x - z| < \delta_2$, then

$$|\tau_2 - \lambda_2| < C e^{(\tau_2 \vee \lambda_2)}(FC(1 + G) + G)|x - z|, \tag{3. 5}$$

$$|x_2 - z_2| < (FC + 1) e^{L(\tau_2 \vee \lambda_2)}(FC(1 + G) + G)|x - z| \tag{3. 6}$$

*Remark2: The Theorem 2 is also found in the literature [1], here we note that the first hitting time $T(x)$ in the proof of Lemma 3 of [1] is different with $t(x, u)$ in this Theorem.*

**(Proof)** without loss of generality let $\tau_1 < \lambda_1$.

Then $t(y_{x_1'}(\lambda_1 - \tau_1, u), u) = \tau_2 - \lambda_1$ and if $|y_{x_1'}(\lambda_1 - \tau_1, u) - z_1'| < \delta_1$ for $\delta_1 > 0$ which is defined as in the Theorem 1, then

$$|\tau_2 - \lambda_2| < C e^{L(\tau_2 - \lambda_1 \vee \lambda_2 - \lambda_1)} |y_{x_1'}(\lambda_1 - \tau_1, u) - z_1'|$$

by using (3. 4).

In the right side of above inequality

$$\left|y_{x_1'}(\lambda_1-\tau_1, u)-z_1'\right| \leq \left|y_{x_1'}(\lambda_1-\tau_1, u)-x_1'\right|+\left|x_1'-z_1'\right|$$

and by using (3. 2), (3. 4) again we have

$$\left|y_{x_1'}(\lambda_1-\tau_1, u)-x_1'\right| \leq F\left|\tau_1-\lambda_1\right| < FCe^{L\lambda_1}|x-z|.$$

Also

$$\left|x_1'-z_1'\right| < G|x_1-z_1| < G(1+FC)e^{L\lambda_1}|x-z|,$$

by using the Lipschitz continuity of $g$ and (3. 3).

Hence combining two inequalities, we get

$$\left|y_{x_1'}(\lambda_1-\tau_1, u)-z_1'\right| < \left(FCe^{L\lambda_1}+G(1+FC)e^{L\lambda_1}\right)|x-z|$$
$$= e^{L\lambda_1}\left(FC(1+G)+G\right)|x-z| \qquad (3.7)$$

Now we define $\delta_2 = \min\{\delta_1, \dfrac{e^{-L\lambda_1}\delta_1}{FC(1+G)+G}\}$.

If $|x-z|<\delta_2$, then $\left|y_{x_1'}(\lambda_1-\tau_1, u)-z_1'\right|<\delta_1$ which implies (3. 5) is proved.

Also

$$|x_2-z_2| = \left|y_{x_1'}(\tau_2-\tau_1, u)-y_{z_1'}(\lambda_2-\lambda_1, u)\right|$$
$$\leq \left|y_{x_1'}(\tau_2-\tau_1, u)-y_{z_1'}(\tau_2-\lambda_1, u)\right|+\left|y_{z_1'}(\tau_2-\lambda_1, u)-y_{z_1'}(\lambda_2-\lambda_1, u)\right|$$
$$= \left|y_{y_{x_1'}(\lambda_1-\tau_1)}(\tau_2-\lambda_1, u)-y_{z_1'}(\tau_2-\lambda_1, u)\right|+\left|y_{z_1'}(\tau_2-\lambda_1, u)-y_{z_1'}(\lambda_2-\lambda_1, u)\right|$$
$$< e^{L(\tau_2-\lambda_1)}\left|y_{x_1'}(\lambda_1-\tau_1, u)-z_1'\right|+F|\tau_2-\lambda_2|$$
$$< e^{L\tau_2}\left(FC(1+G)+G\right)|x-z|+CFe^{L(\tau_2\vee\lambda_2)}\left(FC(1+G)+G\right)|x-z|$$
$$\leq e^{L(\tau_2\vee\lambda_2)}(FC+1)\left(FC(1+G)+G\right)|x-z|,$$

where we used the semigroup property $y_{x_1'}(\tau_2-\tau_1) = y_{y_{x_1'}(\lambda_1-\tau_1)}(\tau_2-\lambda_1, u)$ in the second equality and (3. 1), (3. 2) in the second inequality also (3. 5) and (3. 7) in the third inequality. □

Repeating similar procedure in the sequel discussion we can get a generalized estimation for $i$ th hitting times of trajectories when they hit $A$.

**Theorem 3:** let the $i$ th hitting times of trajectories starting from $x, z$ and evolving with fixed control $u$ be $\tau_i, \lambda_i$, and $x_i, z_i$ be the points where the trajectories hit $A$.

Under the assumptions of Theorem 1, if $|x-z|<\delta_i$, then

$$|\tau_i-\lambda_i| < Ce^{L(\tau_i\vee\lambda_i)}P^{i-1}|x-z| \qquad (3.8)$$
$$|x_i-z_i| < e^{L(\tau_i\vee\lambda_i)}(FC+1)P^{i-1}|x-z|. \qquad (3.9)$$

Here $\delta_i = \min\{\delta_1, \delta_2, \cdots, \dfrac{\delta_1 e^{-L(\tau_i\vee\lambda_i)}}{P^{i-1}}\}$ and $P = FC(1+G)+G$.

The proof is similar with theorem 2.

**Corollary:** Assume ($A_1$)-($A_9$) to be satisfied. Let the reach times of trajectories starting from $x, z$ and evolving with fixed control $u$ be $\overline{\tau} = t_x(u), \overline{\lambda} = t_z(u)$

Denote $\overline{x} = y_x(\overline{\tau}, u), \overline{z} = y_z(\overline{\lambda}, u)$.

Then we get

$$|\bar{\tau}-\bar{\lambda}|<Ce^{L(\bar{\tau}\vee\bar{\lambda})}P^n|x-z| \tag{3.10}$$

$$|\bar{x}-\bar{z}|<e^{L(\bar{\tau}\vee\bar{\lambda})}(FC+1)P^n|x-z| \tag{3.11}$$

for $|x-z|<\bar{\delta}$, where $\bar{\delta}=\min\{\delta_1,\cdots,\delta_n,\dfrac{\delta_1 e^{-L(\bar{\tau}\vee\bar{\lambda})}}{P^n}\}$ and $n$ is the number of times that the trajectories hit the jump set $A$ before they arrive at a target set.

The proof goes through almost the same way as in Theorem 1.

However if we give $\bar{\delta}$ defined above, then $\left|y_{x_n'}(\lambda_n-\tau_n,u)-z_n'\right|<\delta_1$ and

$\left|\bar{x}-y_{z_n'}(\bar{\tau}-\lambda_n,u)\right|<\delta^*<\dfrac{\beta}{2}$ (where $\delta^*$ is defined in Lemma), so $y_{z_n'}(\bar{\tau}-\lambda_n,u)\in B\left(\bar{x},\dfrac{\beta}{2}\right)$.

From the lemma it is clear that if the trajectory starting from $x_n'$ arrives at a target set then the trajectory starting from $z_n'$ also arrives at a target set without hitting $A_j$, which means the number of times of the trajectories starting from $x,z$ respectively and evolving with $u$ is equal if the initial points are closer enough to be satisfied $|x-z|<\bar{\delta}$.

Now we deal with the continuity of the value function through the following theorem.

**Theorem 4**: Under the assumptions ($A_1$)-($A_{12}$), the value function $V$ defined in (2.5) is bounded and if $V$ is lower semicontinuous at each point of $\partial\Gamma$, then it is locally Holder continuous on $\Omega\setminus\bar{\Gamma}$ with respect to the initial point.

**(Proof)** First of all we show the value function is bounded.

From the definition of the value function we get for any $u(\cdot)\in\mathcal{U}$, and $v\in V_1$

$$V(x,q)\leq\int_0^{t_x(u)}K(y_x(t,u),q(t),u(t))e^{-\lambda t}dt+\sum_{\tau_i<t_x(u)}C_a(x_i,q_i,v)e^{-\lambda\tau_i}$$
$$+e^{-\lambda t_x(u)}h(y_x(t_x(u),u)).$$

By assumptions ($A_8$)-($A_{12}$)

$V(x,q)\leq\dfrac{K_0}{\lambda}+C_0\sum_{\tau_i<t_x(u)}e^{-\lambda\tau_i}+H\leq\dfrac{K_0}{\lambda}+C_0e^{-\lambda\tau_1}\dfrac{1}{1-e^{-\frac{\lambda\beta}{F}}}+H$ , this proves the value function is bounded.

Here we used $h(x)$ is bounded by $H$ on $\partial\Gamma$, that is $|h(x)|<H$ $(x\in\partial\Gamma)$

Now we prove the continuity of the value function defined in (2.5).

Let choose $x,z\in\Omega\setminus\Gamma$ to be $|x-z|<\bar{\delta}$.

Considering the definition of the value function we assume that the controller chooses not to make any controlled jump in $C$.

Then by the definition of $V$ for any $\varepsilon>0$, there exist some control $u(\cdot)\in\mathcal{U}$, $v\in V_1$ such that the following inequality holds:

$$V(z,q)\geq\int_0^{t_z(u)}K(y_z(t,u),q(t),u(t))e^{-\lambda t}dt+\sum_{\lambda_i<t_z(u)}C_a(z_i,q_i,v)e^{-\lambda\lambda_i}$$
$$+e^{-\lambda t_z(u)}h(y_z(t_z(u),u))-\varepsilon,$$

where $\{\tau_i\}$, $\{\lambda_i\}$ are sequences of hitting times of the trajectories starting from $x,z$ and evolving with control $u(\cdot)\in\mathcal{U}$, $v\in V_1$

Also

$$V(x, q) \leq \int_0^{t_x(u)} K(y_x(t, u), q(t), u(t))e^{-\lambda t}dt + \sum_{\tau_i < t_x(u)} C_a(x_i, q_i, v)e^{-\lambda \tau_i}$$
$$+ e^{-\lambda t_x(u)} h(y_x(t_x(u), u))$$

Let $t_x(u) = \tilde{s}$, $t_z(u) = \tilde{t}$.

If $\tilde{s} = \tilde{t} = +\infty$ then (2.5) is equal to the infinite horizon hybrid optimal control problem for which the continuity of the value function is proved in [1].

Assume $\tilde{s} \leq \tilde{t} < +\infty$.

Define a new control $\tilde{u}(t) = u(t + \tilde{s})$ and denote the points of trajectories at time $\tilde{s}$ by
$\tilde{z} = y_z(\tilde{s}, u)$, $\tilde{x} = y_x(\tilde{s}, u) \in \partial \Gamma$

Then we have

$$V(z, q) \geq \int_0^{\tilde{s}} K(y_z(t, u), q(t), u(t))e^{-\lambda t}dt + \sum_{\lambda_i < \tilde{s}} C_a(z_i, q_i, v)e^{-\lambda \lambda_i}$$
$$+ \int_{\tilde{s}}^{\tilde{t}} K(y_z(t, u), q(t), u(t))e^{-\lambda t}dt + \sum_{\tilde{s} \leq \lambda_i < \tilde{t}} C_a(z_i, q_i, v)e^{-\lambda \lambda_i} + e^{-\lambda \tilde{t}} h(y_z(\tilde{t}, u)) - \varepsilon$$
$$= \int_0^{\tilde{s}} K(y_z(t, u), q(t), u(t))e^{-\lambda t}dt + \sum_{\lambda_i < \tilde{s}} C_a(z_i, q_i, v)e^{-\lambda \lambda_i} + e^{-\lambda \tilde{s}} J(\tilde{z}, q(\tilde{s}), \tilde{u}, v)$$
$$\geq \int_0^{\tilde{s}} K(y_z(t, u), q(t), u(t))e^{-\lambda t}dt + \sum_{\lambda_i < \tilde{s}} C_a(z_i, q_i, v)e^{-\lambda \lambda_i} + e^{-\lambda \tilde{s}} V(\tilde{z}, q(\tilde{s}))$$

Hence
$$V(x, q) - V(z, q) \leq$$
$$\int_0^{\tilde{s}} |K(y_x(t, u), q(t), u(t)) - K(y_z(t, u), q(t), u(t))| e^{-\lambda t}dt$$
$$+ \sum_{\tau_i, \lambda_i < \tilde{s}} |C_a(z_i, q_i, v) - C_a(x_i, q_i, v)| e^{-\lambda(\tau_i \vee \lambda_i)} \quad (3.12)$$
$$+ e^{-\lambda \tilde{s}} (h(\tilde{x}) - V(\tilde{z}, q(\tilde{s}))) + \varepsilon \quad (\tilde{x} \in \partial \Gamma)$$

Let $I = |K(y_x(t, u), q(t), u(t)) - K(y_z(t, u), q(t), u(t))|$ and consider $\int_0^{\tilde{s}} I e^{-\lambda t} dt$.

If we split the interval $[0, \tilde{s}]$ into a small intervals by $\tau_i$, $\lambda_i$ so that where the system does not make any jumps, then we can calculate $\int_0^{\tilde{s}} I e^{-\lambda t} dt$ separately in three cases.

Case 1: $\int_{\tau_i}^{\lambda_i} I e^{-\lambda t} dt$ or $\int_{\lambda_i}^{\tau_i} I e^{-\lambda t} dt$

In this case by using the Lipschitz continuity of $K$ we have
$$\int_{\tau_i}^{\lambda_i} I e^{-\lambda t} dt \leq 2K_0 |\tau_i - \lambda_i| \quad (3.13)$$

We can go through the same procedure for the calculation of $\int_{\lambda_i}^{\tau_i} I e^{-\lambda t} dt$ that is

$$\int_{\lambda_i}^{\tau_i} I e^{-\lambda t} dt \leq 2K_0 |\tau_i - \lambda_i|$$

Case 2: $\int_{\lambda_i}^{\tau_{i+1}} I e^{-\lambda t} dt$ or $\int_{\lambda_i}^{\lambda_{i+1}} I e^{-\lambda t} dt$

In this case it means that $\tau_i < \lambda_i$ so using the Lipschitz continuity of $K$

$$\int_{\lambda_i}^{\tau_{i+1}} I e^{-\lambda t} dt = \int_{\lambda_i}^{\tau_{i+1}} \left| K\left(y_{x_i'}(t-\tau_i, u), q(t), u(t)\right) - K\left(y_{z_i'}(t-\lambda_i, u), q(t), u(t)\right) \right| e^{-\lambda t} dt$$

$$\leq K_1 \int_{\lambda_i}^{\tau_{i+1}} \left| y_{x_i'}(t-\tau_i, u) - y_{z_i'}(t-\lambda_i, u) \right| e^{-\lambda t} dt$$

In the integral of right side of above inequality

$$\left| y_{x_i'}(t-\tau_i, u) - y_{z_i'}(t-\lambda_i, u) \right| = \left| y_{y_{x_i'}(\lambda_i - \tau_i, u)}(t-\lambda_i, u) - y_{z_i'}(t-\lambda_i, u) \right|$$

$$\leq e^{L(t-\lambda_i)} \left| y_{x_i'}(\lambda_i - \tau_i, u) - z_i' \right|,$$

where we used the semigroup property $y_{x_i'}(t-\tau_i, u) = y_{y_{x_i'}(\lambda_i - \tau_i, u)}(t-\lambda_i, u)$ in the first equality and (3.1) in the next inequality.

On the other hand

$$\left| y_{x_i'}(\lambda_i - \tau_i, u) - z_i' \right| \leq \left| y_{x_i'}(\lambda_i - \tau_i, u) - x_i' \right| + \left| x_i' - z_i' \right|$$

and from (3.2) and (3.8)

$$\left| y_{x_i'}(\lambda_i - \tau_i, u) - x_i' \right| \leq F |\lambda_i - \tau_i| \leq F C e^{L\lambda_i} P^{i-1} |x-z|,$$

Also by using Lipschitz continuity of $g$ and (3.9)

$$\left| x_i' - z_i' \right| \leq G |x_i - z_i| \leq G e^{L\lambda_i} (FC+1) P^{i-1} |x-z|.$$

Combining two inequalities we get

$$\left| y_{x_i'}(\lambda_i - \tau_i, u) - x_i' \right| = P^i e^{L\lambda_i} |x-z|$$

and

$$\left| y_{x_i'}(t-\tau_i, u) - y_{z_i'}(t-\lambda_i, u) \right| = P^i e^{Lt} |x-z|$$

Hence

$$\int_{\lambda_i}^{\tau_{i+1}} I e^{-\lambda t} dt = K_1 P^i |x-z| \int_{\lambda_i}^{\tau_{i+1}} e^{(L-\lambda)t} dt. \tag{3.14}$$

Similarly for $\int_{\lambda_i}^{\lambda_{i+1}} I e^{-\lambda t} dt$ we have

$$\int_{\lambda_i}^{\lambda_{i+1}} I e^{-\lambda t} dt = K_1 P^i |x-z| \int_{\lambda_i}^{\lambda_{i+1}} e^{(L-\lambda)t} dt.$$

Case 3: $\int_{\tau_i}^{\tau_{i+1}} I e^{-\lambda t} dt$ or $\int_{\tau_i}^{\lambda_{i+1}} I e^{-\lambda t} dt$

Likewise in case 2 it can be written

$$\int_{\tau_i}^{\tau_{i+1}} I e^{-\lambda t} dt \leq K_1 \int_{\tau_i}^{\tau_{i+1}} \left| y_{z_1'}(t-\lambda_i, u) - y_{x_i'}(t-\tau_i, u) \right| e^{-\lambda t} dt$$

$$\leq K_1 \int_{\tau_i}^{\tau_{i+1}} \left| y_{y_{z_1'}(\tau_i - \lambda_i)}(t-\tau_i, u) - y_{x_i'}(t-\tau_i, u) \right| e^{-\lambda t} dt$$

$$\leq K_1 \int_{\tau_i}^{\tau_{i+1}} e^{L(t-\tau_i)} P^i e^{L\tau_i} |x-z| e^{-\lambda t} dt \leq K_1 P^i |x-z| \int_{\tau_i}^{\tau_{i+1}} e^{(L-\lambda)t} dt \quad (3.15)$$

Also in the same way we have

$$\int_{\tau_i}^{\lambda_{i+1}} I e^{-\lambda t} dt = K_1 P^i |x-z| \int_{\tau_i}^{\lambda_{i+1}} e^{(L-\lambda)t} dt$$

Now we sum up all the integrals to calculate $\int_0^{\tilde{s}} I e^{-\lambda t} dt$.

Let the trajectories starting from initial points make $n$ times of jumps on $\partial A$ until it reaches a target set.

Define $\tau_0 = 0$, $\tau_{n+1} = \tilde{s}$, then using (3.8), (3.13), (3.14), (3.15)

$$\int_o^{\tilde{s}} I e^{-\lambda t} dt \leq \sum_{i=1}^n 2K_0 C P^{i-1} e^{L\tilde{s}} |x-z| + \sum_{i=0}^n K_1 P^i |x-z| \int_o^{\tilde{s}} e^{(L-\lambda)t} dt$$

$$= \begin{cases} 2K_0 C e^{L\tilde{s}} \left( \dfrac{P^n - 1}{P-1} \right) |x-z| + K_1 \left( \dfrac{P^{n+1}-1}{P-1} \right) \cdot \dfrac{e^{(L-\lambda)\tilde{s}} - 1}{L-\lambda} |x-z| & L \neq \lambda \\ 2K_0 C e^{L\tilde{s}} \left( \dfrac{P^n - 1}{P-1} \right) |x-z| + K_1 \left( \dfrac{P^{n+1}-1}{P-1} \right) \cdot \tilde{s} |x-z| & L = \lambda \end{cases} \quad (3.16)$$

On the other hand by Lipschitz continuity of $C_a$ and (3.9) we get

$$\sum_{i=1}^n \left| C_a(x_i, q_i, v) - C_a(z_i, q_i, v) \right| e^{-\lambda(\tau_i \vee \lambda_i)} \leq \sum_{i=1}^n 2C_1 |x_i - z_i| e^{-\lambda(\tau_i \vee \lambda_i)}$$

$$\leq 2C_1 \sum_{i=1}^n (FC+1) e^{L\tilde{s}} P^{i-1} |x-z| \leq 2C_1 (FC+1) e^{L\tilde{s}} \left( \frac{P^n - 1}{P-1} \right) |x-z|$$

We can think some constant $a$ so that the following inequality holds

$$\frac{P^n - 1}{P - 1} < aP^n.$$

In fact it is possible if we give $F$, $G$, $C$ large constant then $P$ is also large enough to be satisfied above inequality. (For example it is sufficient when we give $a = \dfrac{1}{P-1}$.)

Also we have $\tau_{i+1} - \tau_i \geq \dfrac{\beta}{F}$, $\tilde{s} \geq \tau_{n+1} - \tau_1 \geq n \cdot \dfrac{\beta}{F}$, $n \leq \dfrac{\tilde{s}F}{\beta}$, then by (3.12), (3.16) it follows that

$$V(x, q) - V(z, q) \leq \begin{cases} 2aK_0 Ce^{L\tilde{s}} P^{\frac{\tilde{s}F}{\beta}} |x-z| + aK_1 P \cdot \dfrac{e^{L\tilde{s}} P^{\frac{\tilde{s}F}{\beta}}}{L-\lambda} |x-z| + \\ \quad e^{-\lambda \tilde{s}} \left( h(\tilde{x}) - V(\tilde{z}, q(\tilde{s})) \right) + \varepsilon & L \neq \lambda \\ 2aK_0 Ce^{L\tilde{s}} P^{\frac{\tilde{s}F}{\beta}} |x-z| + aK_1 P \cdot P^{\frac{\tilde{s}F}{\beta}} \cdot \tilde{s} |x-z| \\ \quad + e^{-\lambda \tilde{s}} \left( h(\tilde{x}) - V(\tilde{z}, q(\tilde{s})) \right) + \varepsilon & L = \lambda \end{cases} \quad (3.17)$$

Now we fix some $\theta \in (0, 1)$ and define $\delta = \min\left\{ \bar{\delta}, \left( e^{(L+1)\tilde{s}} \cdot P^{\frac{\tilde{s}F}{\beta}} \right)^{-\frac{1}{\theta}} \right\}$.

Let $T = \dfrac{-\theta \ln \delta}{L + 1 + \dfrac{F}{\beta} \ln P}$.

Then for any initial points $x, z$ ($|x-z| < \delta$),

$$T = \frac{-\theta \ln \delta}{L + 1 + \dfrac{F}{\beta} \ln P} < \frac{-\theta \ln \delta}{L + \dfrac{F}{\beta} \ln P} < \frac{-\theta \ln |x-z|}{L + \dfrac{F}{\beta} \ln P}$$

and from the fact

$$e^{LT} P^{\frac{TF}{\beta}} < |x-z|^{-\theta}$$

$$\Leftrightarrow LT + \frac{TF}{\beta} \ln P < -\theta \ln |x-z|$$

$$\Leftrightarrow T < \frac{-\theta \ln |x-z|}{L + \dfrac{F}{\beta} \ln P}$$

we have $e^{LT} P^{\frac{TF}{\beta}} < |x-z|^{-\theta}$. (Here we assume $L + \dfrac{F}{\beta} \ln P > 0$)

Also

$$\tilde{s} < \frac{-\theta \ln \delta}{L + 1 + \dfrac{F}{\beta} \ln p}$$

$$\Leftrightarrow \tilde{s}(L+1) + \tilde{s} \frac{F}{\beta} \ln P < -\theta \ln \delta$$

$$\Leftrightarrow e^{(L+1)\tilde{s}} P^{\frac{\tilde{s}F}{\beta}} < \delta^{-\theta} \Leftrightarrow \delta < \left( e^{(L+1)\tilde{s}} \cdot P^{\frac{\tilde{s}F}{\beta}} \right)^{-\frac{1}{\theta}}$$

then $\tilde{s} < T$ by the definition of $\delta$.

Therefore (2.17) can be written as follows

$$V(x, q) - V(z, q) \leq \begin{cases} M_1 |x-z|^{1-\theta} + e^{-\lambda \tilde{s}} \left( h(\tilde{x}) - V(\tilde{z}, q(\tilde{s})) \right) + \varepsilon & L \neq \lambda \\ M_2 |x-z|^{1-\theta} + e^{-\lambda \tilde{s}} \left( h(\tilde{x}) - V(\tilde{z}, q(\tilde{s})) \right) + \varepsilon & L = \lambda \end{cases} \quad (3.18)$$

Here $M_1, M_2$ are some constants.

To prove the local Holder continuity of the value function we observe that it is continuous on $\partial \Gamma$.

We can choose such small neighborhood of a target set that no jumps can make in any points of the neighborhood.

Then for control $u$

$$V(\tilde{z}, q(\tilde{s})) - h(\tilde{x}) \le \int_0^{t_{\tilde{z}}(u)} K(y_{\tilde{z}}(t, u), q(t), u(t)) e^{-\lambda t} dt + h(y_{\tilde{z}}(t_{\tilde{z}}(u), u)) - h(\tilde{x})$$

We show that the first term in the right side of above inequality tends to zero as $\tilde{z} \to \tilde{x}$.

Define $T(x) = \inf_{u \in \mathcal{U}} t_x(u)$ then we pick $u^* \in \mathcal{U}$ such that $t_{\tilde{z}}(u) \le T(\tilde{z}) + \varepsilon$ for any $\varepsilon > 0$.

Then

$$\int_0^{t_{\tilde{z}}(u)} K(y_{\tilde{z}}(t, u), q(t), u(t)) e^{-\lambda t} dt \le K_0 t_{\tilde{z}}(u) \le K_0 (T(\tilde{z}) + \varepsilon),$$

where $K_0$ is defined in assumption ($A_{11}$)

The right side tends to zero as $\tilde{z} \to \tilde{x}$ and $\varepsilon \to 0^+$ thanks the continuity of $T$ (see the proposition 1. 2 in chapter IV of [3]).

Also $|y_{\tilde{z}}(t_{\tilde{z}}(u), u) - \tilde{x}| \le |y_{\tilde{z}}(t_{\tilde{z}}(u), u) - \tilde{z}| + |\tilde{z} - \tilde{x}| \to 0$ and this means

$$h(y_{\tilde{z}}(t_{\tilde{z}}(u), u)) - h(\tilde{x}) \to 0$$

by the continuity of $h$ on $\Gamma$.

As a result, we get $\limsup_{\tilde{z} \to \tilde{x}} V(\tilde{z}, q(\tilde{s})) \le h(\tilde{x})$, $\tilde{x} \in \partial \Gamma$.

Thanks to the assumption that $V$ is lower semicontinuous on $\partial \Gamma$, it is continuous on $\partial \Gamma$.

Therefore (3. 18) tends to zero as $x \to z$, because $|\tilde{x} - \tilde{z}| \le e^{\tilde{s}} |x - z| \to 0$ as $x \to z$ and by the continuity of the value function on $\partial \Gamma$ it follows that $h(\tilde{x}) - V(\tilde{z}, q(\tilde{s})) \to 0$ as $x \to z$.

Now consider in the case $\tilde{s} > \tilde{t}$.

In this case from the dynamic programming theory (we prove it in the next section) we have

$$V(x, q) \le \int_0^{\tilde{t}} K(y_x(t, u), q(t), u(t)) e^{-\lambda t} dt + \sum_{\tau_i < \tilde{t}} C_a(x_i, q_i, v) e^{-\lambda \tau_i} + e^{-\lambda \tilde{t}} V(y_x(\tilde{t}, u), q(\tilde{t}))$$

$$\therefore V(x, q) - V(z, q) \le \int_0^{\tilde{t}} |K(y_x(t, u), q(t), u(t)) - K(y_z(t, u), q(t), u(t))| e^{-\lambda t} dt$$

$$+ \sum_{\tau_i, \lambda_i < \tilde{s}} |C_a(z_i, q_i, v) - C_a(x_i, q_i, v)| e^{-\lambda(\tau_i \vee \lambda_i)}$$

$$+ e^{-\lambda \tilde{t}} |V(y_x(\tilde{t}, u), q(\tilde{t})) - h(y_z(\tilde{t}, u))| + \varepsilon \quad (y_z(\tilde{t}, u) \in \partial \Gamma)$$

Taking throughout the same procedure as in the case $\tilde{s} \le \tilde{t}$, we can prove the continuity of the value function on $\Omega$ under the assumption that it is lower semicontinuous on $\partial \Gamma$. □

## 4. Dynamic Programming Principle and Hamilton-Jacobi variational inequality

First we show the following theorem which is called for Dynamic Programming Principle.
**Theorem 5:** Let $V(x, q)$ be the value function defined in (2. 5).

1) If $\tau_1$ is the first hitting time that trajectory starting from the initial point $x$ and evolving with control $u$ hits $A$, then

$$V(x, q) = \inf_u \{ \int_0^{\tau_1} K(y_x(t, u), q(t), u(t)) e^{-\lambda t} dt + e^{-\lambda \tau_1} MV(x_1, q_1) \}, \quad (4. 1)$$

where $MV(x, q) = \inf_{v \in V_1}\{V(g(x, q, v), q) + C_a(x, q, v)\}$.

2) If $\xi_1$ is the first time that the controller chooses to make a jump in $C$, then

$$V(x, q) = \inf_u \{\int_0^{\xi_1} K(y_x(t, u), q(t), u(t))e^{-\lambda t}dt + e^{-\lambda \xi_1}NV(y(\xi_1), q(\xi_1))\}, \quad (4.2)$$

where $NV(x, q) = \inf_{(x', q') \in D \times I}\{V(x', q') + C_c(x, q, x', q')\}$.

3) For any $T > 0$

$$V(x) = \inf_{u, v, \xi_i, y(\xi_i)'} \{\int_0^{T \wedge t_x(u)} K(y_x(t, u), q(t), u(t))e^{-\lambda t}dt + \sum_{\tau_i < T} C_a(x_i, q_i, v)e^{-\lambda \tau_i}$$

$$+ \sum_{\xi_i < T} C_c\left(y(\xi_i), q(\xi_i), y(\xi_i)', q(\xi_i)'\right)e^{-\lambda \xi_i}$$

$$+ e^{-\lambda(T \wedge t_x(u))}V(y_x(T \wedge t_x(u), u), q(T \wedge t_x(u)))\}, \quad (4.3)$$

where $a \wedge b = \min\{a, b\}$.

**(Proof)** First let's prove 1).

If $\tau_1$ is the first hitting time of $A$, then

$$V(x, q) \leq \int_0^{\tau_1} K(y_x(t, u), q(t), u(t))e^{-\lambda t}dt + C_a(x_1, q_1, v)e^{-\lambda \tau_1}$$

$$+ (\int_{\tau_1}^{t_x(u)} K(y_x(t, u), q(t), u(t))e^{-\lambda t}dt + \sum_{\tau_1 < \tau_i < t_x(u)} C_a(x_i, q_i, v)e^{-\lambda \tau_i}$$

$$+ \sum_{\xi_i < t_x(u)} C_c\left(y(\xi_i), q(\xi_i), y(\xi_i)', q(\xi_i)'\right)e^{-\lambda \xi_i} + e^{-\lambda t_x(u)}h(y_x(t_x(u), u)))$$

From the third term in the right side of above inequality we change the variable $t' = t - \tau_1$ in the bracket and after that taking infimum in the bracket over the control variables we have

$$V(x, q) \leq \int_0^{\tau_1} K(y_x(t, u), q(t), u(t))e^{-\lambda t}dt + e^{-\lambda \tau_1}C_a(x_1, q_1, v)$$

$$+ e^{-\lambda \tau_1}V(g(x_1, q_1, v), q(\tau_1))$$

Here again taking the infimum in the last two terms over the discrete control $v \in V_1$ and taking infimum over the control $u \in \mathcal{U}$ again, then we get

$$V(x, q) \leq \inf_u \{\int_0^{\tau_1} K(y_x(t, u), q(t), u(t))e^{-\lambda t}dt + e^{-\lambda \tau_1}MV(x_1, q_1)\}$$

For the reverse inequality if we give any $\varepsilon > 0$, then thanks to the definition of the value function there exist some control $\bar{\theta} = \left(\bar{u}, \bar{v}, \bar{\xi}_i, y(\bar{\xi}_i)'\right)$ such that

$$V(x, q) \geq \int_0^{\bar{\tau}_1} K(y_x(t, \bar{u}), q(t), \bar{u}(t))e^{-\lambda t}dt + C_a\left(\bar{x}_1, \bar{q}_1, \bar{v}\right)e^{-\lambda \bar{\tau}_1} +$$

$$(\int_{\bar{\tau}_1}^{t_x(\bar{u})} K(y_x(t, \bar{u}), q(t), \bar{u}(t))e^{-\lambda t}dt + \sum_{\bar{\tau}_1 < \bar{\tau}_i < t_x(\bar{u})} C_a\left(\bar{x}_i, \bar{q}_i, \bar{v}\right)e^{-\lambda \bar{\tau}_i}$$

$$+ \sum_{\bar{\xi}_i < t_x(\bar{u})} C_c\left(y(\bar{\xi}_i), q(\bar{\xi}_i), y(\bar{\xi}_i)', q(\bar{\xi}_i)'\right)e^{-\lambda \bar{\xi}_i} + e^{-\lambda t_x(\bar{u})}h(y_x(t_x(\bar{u}), \bar{u}))) - \varepsilon$$

We change the variable $t' = t - \tau_1$ in the bracket and taking infimum over the control variables in the bracket, then we get

$$V(x, q) \geq \int_0^{\bar{\tau}_1} K(y_x(t, \bar{u}), q(t), \bar{u}(t)) e^{-\lambda t} dt + e^{-\lambda \bar{\tau}_1} C_a(\bar{x}_1, \bar{q}_1, \bar{v})$$

$$+ e^{-\lambda \bar{\tau}_1} V(g(\bar{x}_1, \bar{q}_1, \bar{v}), q(\bar{\tau}_1)) - \varepsilon$$

Taking infimum over the discrete control $v \in V_1$ in the last two terms and taking infimum again over the control $u \in \mathcal{U}$, then

$$V(x) \geq \inf_u \{\int_0^{\tau_1} K(y_x(t, u), q(t), u(t)) e^{-\lambda t} dt + e^{-\lambda \tau_1} MV(x_1, q_1)\} + \varepsilon$$

Hence as $\varepsilon \to 0$ we get the reverse inequality.

We can prove 2) in a similar way as in 1).

Now proceed to prove 3).

If $T \geq t_x(u)$ then (4. 3) is equivalent to the definition of the value function (2. 5).

Let $T < t_x(u)$. (4. 3) can be written as follows.

$$V(x, q) = \inf_{u, v, \xi_i, y(\xi_i)'} \{\int_0^T K(y_x(t, u), q(t), u(t)) e^{-\lambda t} dt + \sum_{\tau_i < T} C_a(x_i, q_i, v) e^{-\lambda \tau_i}$$

$$+ \sum_{\xi_i < T} C_c(y(\xi_i), q(\xi_i), y(\xi_i)', q(\xi_i)') e^{-\lambda \xi_i} + e^{-\lambda(T)} V(y_x(T, u), q(T))\}$$

By the definition of the value function for any control $u$ we have

$$V(x, q) \leq \int_0^T K(y_x(t, u), q(t), u(t)) e^{-\lambda t} dt + \sum_{\tau_i < T} C_a(x_i, q_i, v) e^{-\lambda \tau_i}$$

$$+ \sum_{\xi_i < T} C_c(y(\xi_i), q(\xi_i), y(\xi_i)', q(\xi_i)') e^{-\lambda \xi_i}$$

$$+ (\int_T^{t_x(u)} K(y_x(t, u), q(t), u(t)) e^{-\lambda t} dt + \sum_{T \leq \tau_i < t_x(u)} C_a(x_i, q_i, v) e^{-\lambda \tau_i}$$

$$+ \sum_{T \leq \xi_i < t_x(u)} C_c(y(\xi_i), q(\xi_i), y(\xi_i)', q(\xi_i)') e^{-\lambda \xi_i} + e^{-\lambda t_x(u)} h(y_x(t_x(u), u)))$$

We change the variables $t' = t - T$ in the bracket in the right side of above inequality and taking infimum over the control variables. Hence

$$V(x, q) \leq \int_0^T K(y_x(t, u), q(t), u(t)) e^{-\lambda t} dt + \sum_{\tau_i < T} C_a(x_i, q_i, v) e^{-\lambda \tau_i}$$

$$+ \sum_{\xi_i < T} C_c(y(\xi_i), q(\xi_i), y(\xi_i)', q(\xi_i)') e^{-\lambda \xi_i} + e^{-\lambda T} V(y_x(T, u), q(T))$$

Continue taking infimum over the control variables in the right side of above inequality to get

$$V(x, q) \leq \inf_{u, v, \xi_i, y(\xi_i)'} \{\int_0^T K(y_x(t, u), q(t), u(t)) e^{-\lambda t} dt + \sum_{\tau_i < T} C_a(x_i, q_i, v) e^{-\lambda \tau_i}$$

$$+ \sum_{\xi_i < T} C_c(y(\xi_i), q(\xi_i), y(\xi_i)', q(\xi_i)') e^{-\lambda \xi_i} + e^{-\lambda T} V(y_x(T, u), q(T))\} \quad (4.4)$$

To get a reverse inequality let give any $\varepsilon > 0$.

By the definition of the value function there exist some control $\bar{\theta} = \left(\bar{u}, \bar{v}, \overline{\xi_i}, y(\overline{\xi_i})'\right)$ such that

$$V(x, q) \geq \int_0^T K(y_x(t, \bar{u}), q(t), \bar{u}(t)) e^{-\lambda t} dt + \sum_{\tau_i < T} C_a(\overline{x_i}, \overline{q_i}, \bar{v}) e^{-\lambda \overline{\tau_i}}$$

$$+ \sum_{\overline{\xi_i} < T} C_c\left(y(\overline{\xi_i}), q(\overline{\xi_i}), y(\overline{\xi_i})', q(\overline{\xi_i})'\right) e^{-\lambda \overline{\xi_i}} + (\int_T^{t_x(\bar{u})} K(y_x(t, \bar{u}), q(t), \bar{u}(t)) e^{-\lambda t} dt$$

$$+ \sum_{T \leq \overline{\tau_i} < t_x(\bar{u})} C_a(\overline{x_i}, \overline{q_i}, \bar{v}) e^{-\lambda \overline{\tau_i}} + \sum_{T \leq \overline{\xi_i} < t_x(\bar{u})} C_c\left(y(\overline{\xi_i}), q(\overline{\xi_i}), y(\overline{\xi_i})', q(\overline{\xi_i})'\right) e^{-\lambda \overline{\xi_i}}$$

$$+ e^{-\lambda t_x(\bar{u})} h(y_x(t_x(\bar{u}), \bar{u}))) - \varepsilon$$

Changing the variable $t' = t - T$ in the bracket from the third term to last one and taking infimum over the control variables, then we have

$$V(x, q) \geq \int_0^T K(y_x(t, \bar{u}), q(t), \bar{u}(t)) e^{-\lambda t} dt + \sum_{\tau_i < T} C_a(\overline{x_i}, \overline{q_i}, v) e^{-\lambda \overline{\tau_i}}$$

$$+ \sum_{\overline{\xi_i} < T} C_c\left(y(\overline{\xi_i}), q(\overline{\xi_i}), y(\overline{\xi_i})', q(\overline{\xi_i})'\right) e^{-\lambda \overline{\xi_i}} + e^{-\lambda T} V(y_x(T, \bar{u}), q(T)) - \varepsilon$$

Taking infimum over the control variables again we get the reverse inequality as $\varepsilon \to 0$

$$V(x, q) \geq \inf_{u, v, \xi_i, y(\xi_i)'} \{\int_0^T K(y_x(t, u), q(t), u(t)) e^{-\lambda t} dt + \sum_{\tau_i < T} C_a(x_i, q_i, v) e^{-\lambda \tau_i}$$

$$+ \sum_{\xi_i < T} C_c\left(y(\xi_i), q(\xi_i), y(\xi_i)', q(\xi_i)'\right) e^{-\lambda \xi_i} + e^{-\lambda T} V(y_x(T, u), q(T))\} \qquad (4.5)$$

From (4.4) and (4.5) we prove 3). □

Now we derive the Hamilton-Jacobi inequality that the value function must be satisfied by using Dynamic Programming Principle.

The value function of the classical reach time problems which is controlled by only continuous control usually is given as a viscosity solution of Dirichlet boundary problem with the terminal cost $h$ as a boundary data.

In hybrid optimal control problem involving both of the continuous and discrete controls the partial differential equation satisfied by the value function is represented as a quasi-variational inequality. Also the Hamiltonian is given as follows

$$H(x, q, p) := \sup_{u \in U} \left\{ \frac{-K(x, q, u) - f(x, q, u) \cdot p}{\lambda} \right\}$$

We prove the following theorem which gives a method to determine the value function as a viscosity solution of the Hamilton-Jacobi variational inequality.

**Theorem 6:** Assume ($A_1$)-($A_{12}$). If the value function of the hybrid optimal control problem with reach time to a target set defined in (2.5) is lower semicontinuous on $\partial \Gamma$, then it satisfies the following quasi-variational inequality in the viscosity sense.

$$\begin{cases} V(x, q) - MV(x, q) = 0 & (x, q) \in A \times I \\ \max\{V(x, q) - NV(x, q), V(x, q) + H(x, q, D_x V(x, q))\} = 0 & (x, q) \in C \times I \\ V(x, q) + H(x, q, D_x V(x, q)) = 0 & (x, q) \in \Omega \setminus (A \cup C \cup \Gamma) \times I \\ V(x, q) - h(x) = 0 & (x, q) \in \partial \Gamma \times \{j\} \end{cases} \qquad (4.6)$$

*Remark3. By subsolution (respectively, supersolution) of Dirichlet condition $V(x, q) = h(x, q)$ on $(x, q) \in \partial\Gamma \times \{j\}$ we mean a function $\leq h$ (respectively, $\geq h$) at each point of $(x, q) \in \partial\Gamma \times \{j\}$*

The proof of the Theorem 6 is essentially the same as the proof of the Theorem 4. 2 of [1] and we omit it.

## 5. Uniqueness

In this section we prove the comparison between any two viscosity solutions to demonstrate the uniqueness of the viscosity solutions of QVI.

The following theorem is the extension of theorem 5. 1 of [1] to a comparison result for a Dirichlet boundary value problem.

**Theorem 7**: Assume ($A_1$)-($A_{12}$). If the value function defined in (2. 5) is continuous then it is the unique solution of (4. 6) in viscosity sense.

*Remark4. The condition of continuity of the value function can be weakened in the way that it is lower semicontinuous on $\partial\Gamma$. In this case V is the unique solution of QVI in the viscosity sense thanks to the Theorem4.*

(**Proof**) The idea of proof is quite similar with the one noted in the proof of Theorem 5. 1 of [1] so we briefly note here comparing with the proof of above Theorem.

Let $u_1, u_2$, bounded and continuous functions on $\overline{\Omega \setminus \Gamma}$, be two viscosity solutions of QVI (given by (4. 6))

Define the auxiliary function $\Phi$ on $\Omega \times \Omega$ as follows:

$$\Phi^q(x, y) = u_1(x) - u_2(y) - \frac{1}{\varepsilon}|x-y|^2 - k(|x|^2 + |y|^2), \quad (x, y) \in \Omega_q \times \Omega_q$$

where $\varepsilon$ and $k$ are small positive parameters to be chosen conveniently.

Now we prove $\sup_q \sup_{\Omega_q \times \Omega_q} \Phi^q(x, y) \leq 0$

Suppose that $\sup_q \sup_{\Omega_q \times \Omega_q} \Phi^q(x, y) = C > 0$ and fix $k$ such that $k < \min\left\{\frac{C}{2}, \frac{C'}{2}\right\}$ where $C'$ is the constant to be chosen such that $C_a, C_c$ are bounded below. (See assumption ($A_{12}$).)

Then we can choose some point $(x_k, y_k) \in \Omega_1 \times \Omega_1$ without loss of generality such that

$\Phi^1(x_k, y_k) > C - k > \frac{C}{2}$ by the definition of supremum.

Let $\Phi^1$ attain its supremum at some points $(x_0, y_0)$ in $\Omega_1 \times \Omega_1$. We can choose this point because thanks to the structure of $\Phi$ and boundedness of $u_1, u_2$, $\Phi^1(x, y) \to -\infty$ as $|x|, |y| \to \infty$.

We first summarize the following estimates needed in the proof of the uniqueness.

**Lemma2**. The following statements are true.

(1) $\frac{|x_0 - y_0|^2}{\varepsilon} \leq D$ for some $D$ independent of $\varepsilon$ and $k$.

(2) $\sqrt{k}|x_0|, \sqrt{k}|y_0| \leq \hat{D}$ for some $\hat{D}$ independent of $\varepsilon$ and $k$.

(3) $\frac{|x_0 - y_0|^2}{\varepsilon} \leq \omega_k(\sqrt{C\varepsilon})$ where $\omega_k$ is the local modulus of continuity of both $u_1, u_2$ in the ball of radius $R(k)$, dependent on $k$, but independent of $\varepsilon$.

The proof is listed in [1].

If $j \neq 1$ that is $\Gamma \not\subset \Omega_1$ then all the statement in the proof of Theorem 5. 1 of [1] can be dropped with no changes in the proof. If not we have some differences.

The only difference in the proof is one of the $(x_0, y_0)$ may belong in $\Gamma$.

If $(x_0, y_0) \in \Gamma \times A$ or $A \times \Gamma$ it is straightforward that it is contradicted to (2) of Lemma2 because of the assumption $(A_8)$.

In the case of $(x_0, y_0) \in \Gamma \times C$ or $C \times \Gamma$ we can get a contradiction in the same way.

Now consider in the case $(x_0, y_0) \in \Gamma \times \Gamma$ or $(x_0, y_0) \in \Gamma \times \Omega \setminus (A \cup C \cup \Gamma)$ or $(x_0, y_0) \in \Omega \setminus (A \cup C \cup \Gamma) \times \Gamma$.

Without loss of generality let $x_0 \in \partial \Gamma$.

$u_1$ and $u_2$ are respectively the viscosity subsolution and supersolution, so we get $u_1 \leq h \leq u_2$ on $\partial \Gamma$. Then

$$\Phi^1(x_0, y_0) = u_1(x_0) - u_2(y_0) - \frac{1}{\varepsilon}|x_0 - y_0|^2 - k\left(|x_0|^2 + |y_0|^2\right)$$
$$\leq u_1(x_0) - u_2(x_0) + u_2(x_0) - u_2(y_0) \leq u_2(x_0) - u_2(y_0) \leq \overline{\omega}\left(\sqrt{C\varepsilon}\right)$$

where $\overline{\omega}$ is the modulus of continuity of $u_2$ in the ball of radius $R(k)$, dependent on $k$, but independent of $\varepsilon$.

By choosing $\varepsilon$ small enough we get a contradiction to $\sup_q \sup_{\Omega_q \times \Omega_q} \Phi^q(x, y) = C > 0$.

If $x_0 \in \partial \Gamma$ we add and subtract $u_1(y_0)$ and use the modulus of continuity of $u_1$ to conclude in the same way.

The other part of the proof is the same as in the proof of Theorem 5.1 of [1]. □